\documentclass[12pt]{article}
\usepackage{graphicx,amssymb}
\usepackage{fullpage}

\newcommand{\text}[1]{\mbox{#1}}
\newcommand{\Expect}{\mathop{\bf E{}}}
\usepackage{amsmath}

\newcommand{\QED}{~~\rule[-1pt]{6pt}{6pt}}
\newcommand{\reals}{{\mbox{\bf R}}}

\newcommand{\symm}{{\mbox{\bf S}}}  

\newcommand{\BEAS}{\begin{eqnarray*}}
\newcommand{\EEAS}{\end{eqnarray*}}
\newcommand{\BEA}{\begin{eqnarray}}
\newcommand{\EEA}{\end{eqnarray}}
\newcommand{\BEQ}{\begin{equation}}
\newcommand{\EEQ}{\end{equation}}
\newcommand{\BIT}{\begin{itemize}}
\newcommand{\EIT}{\end{itemize}}
\newcommand{\BNUM}{\begin{enumerate}}
\newcommand{\ENUM}{\end{enumerate}}

\newcommand{\BA}{\begin{array}}
\newcommand{\EA}{\end{array}}

\newtheorem{theorem}{Theorem}

\newtheorem{corollary}[theorem]{Corollary}

\newtheorem{definition}[theorem]{Definition}

\newtheorem{lemma}[theorem]{Lemma}

\newtheorem{proposition}[theorem]{Proposition}

\newenvironment{proof}{\textbf{Proof.}}{\QED\bigskip}

\begin{document}

\title{A harmonic analysis solution to the static basket arbitrage problem}
\author{Alexandre d'Aspremont \thanks{
ISL \& MS\&E dept, Stanford University, Stanford, CA 94305-4023,
USA. alexandre.daspremont@m4x.org }} \maketitle

\begin{abstract}
We consider the problem of computing upper and lower bounds on the
price of a European basket call option, given prices on other
similar baskets. We focus here on an interpretation of this
program as a generalized moment problem. Recent results by Berg \&
Maserick (1984), Putinar \& Vasilescu (1999) and Lasserre (2001)
on harmonic analysis on semigroups, the $\mathbb{K}$-moment
problem and its applications to optimization, allow us to derive
tractable necessary and sufficient conditions for the absence of
static arbitrage between basket straddles, hence between basket
calls and puts.

\textbf{Keywords}: Semidefinite Programming, Static Arbitrage,
K-Moment Problem, Basket Options.
\end{abstract}

\section{Introduction}

We let $p \in \reals_+^{n+m}$, $K \in \reals^{n+m+1}$, $w_i \in
\reals^n$, $i=0, \ldots, n+m$ and we consider the problem of
computing upper and lower bounds on the price of an European
basket call option with strike $K_0$ and weight vector $w_0$:
\BEQ
\BA{ll}
\mbox{maximize/minimize} & {p_0:=\Expect}_{\nu }(w_0^{T}x-K_0)_+ \\
\mbox{subject to} & {\Expect}_{\nu }(w_i^{T}x-K_i)_+=p_i, \quad
i=1,\ldots,n+m, \label{eq:constraints}
\EA
\EEQ
with respect to all probability measures $\nu$ on the asset price
vector $x \in \reals_+^n$, consistent with the (given) set of
observed prices $p_i$ of options on other baskets.

We implicitly assume that all the options have the same maturity,
and that, without loss of generality, the risk-free interest rate
is zero (we compare prices in the forward market). We seek
non-parametric bounds, i.e., we do not assume any specific model
for the underlying asset prices, our only assumption is the
absence of a static arbitrage today (i.e. the absence of an
arbitrage that only requires trading today and at maturity).

Here, we interpret (\ref{eq:constraints}) as a \textit{generalized
moment problem}. This approach was successfully used in
\cite{Bert00} to get tractable bounds in dimension one and to show
the NP-hardness of the multivariate problem
(\ref{eq:constraints}). NP-hardness means that we have no chance
of finding a direct and efficient method for detecting all
arbitrage opportunities, here instead we look for a sequence of
successively tighter price bounds. This means that outlandish
arbitrage opportunities can be detected at little numerical cost
while detecting finer price discrepancies has a higher theoretical
complexity.

Recent results on multivariate moment problems (see \cite{Schm91},
\cite{Puti99} or \cite{Curt00}), semidefinite programming (see
\cite{Nest94}, \cite{Vand96} and \cite{Nest00}) and harmonic
analysis on semigroups (see \cite{Berg84b} and \cite{Roma03})
allow us to derive static arbitrage price bounds on a set of
products linked by a semigroup structure. The resulting
constraints can be formulated as successively tighter linear
matrix inequalities, hence we can compute increasingly sharp
bounds on the solution to problem (\ref{eq:constraints}) as
solutions of increasingly large semidefinite programs (linear
programs on the cone of positive semidefinite matrices).
Semidefinite programming has been the object of intensive research
since the seminal work of \cite{Nest94} and several numerical
packages (see for example SEDUMI by \cite{stur99}) are now
available to solve these problems very efficiently.

The core of our argument is to substitute to the classical duality
between the cones of probability measures and positive portfolios,
the conic duality between positive definite functions on one hand
and sums of squares on the other. These last two cones have the
advantage of being numerically tractable and lead to exploitable
formulations of the static portfolio super/sub-replication
problems.

A lot of work has been focused on arbitrage bounds in a dynamic
setting, see \cite{ElKa91}, \cite{Karou95}, \cite{Avel95}, and
\cite {Kara98}, among others. Work on the unidimensional static
problem dates back at least to \cite{Bree78} (see also
\cite{Laur00}), both using the positivity of butterfly spread
prices to preclude arbitrage. \cite{Bert00} studied these bounds
together with second order moment constraints and proved the
NP-Hardness of the multivariate problem (\ref{eq:constraints}).
Finally, in a previous paper \cite{dasp03b}, we focused on the
interpretations of problem (\ref{eq:constraints}) as an
\textit{integral transform inversion problem} or a \textit{linear
semi-infinite program}, i.e. a linear program with a finite number
of linear constraints on an infinite dimensional variable, and
used the related theories to compute closed-form solutions for
some particular cases and a linear programming relaxation for the
general case.

The paper is organized as follows. In section two, we describe the
static market structure and start with a brief introduction on
harmonic analysis on semigroups. Based on these results, we then
derive necessary and sufficient conditions for the absence of
arbitrage in the static market, formulated as semidefinite
programs. Finally, in section three, we describe the conic duality
between positive definite functions and sums of squares and use it
to show how a super/sub-replicating portfolio can be constructed
from the solution to the programs of the preceding section.

\section{Static arbitrage constraints}

\subsection{Market structure}
\label{ss-market-structure} We work in a one period framework and
suppose that the market is composed of cash and $n$ underlying
assets $x_i$ for $i=1,\ldots,n$ with $x \in \reals_+^n$. We
suppose that the forward prices of the assets are known and given
by $p_i$, for $i=1,\ldots,n$, hence $w_i$ is the Euclidean basis
and $K_i=0$ for $i=1,\ldots,n$. In addition to these basic
products, there are $m+1$ basket \textit{straddles} on the assets
$x$, with payoff given by $|w_{n+i}^{T}x-K_{n+i}|$, $i=1, \ldots,
m$. Because a straddle is obtained as the sum of a call and a put,
we get the market price of straddles from those of basket calls
and forward contracts by call-put parity.

We will note these payoff functions $e_i$, for $i=0,\ldots,m+n$,
with $e_i(x)=x_i$ for $i=1,\ldots,n$ and
$e_{(n+j)}(x)=|w_i^{T}x-K_i|$ for $j=0,\ldots,m$. In what follows,
we will focus on the Abelian (commutative) semigroup
$(\mathbb{S},\cdot)$ generated by the payoffs $e_i(x)$ for
$i=0,\ldots,m+n$, the cash $1_{\mathbb{S}}$ and their products.

In this one period setting, we will look for conditions that
guarantee the absence of \textit{static arbitrages}, i.e.
arbitrage opportunities that only involve trading today and at
maturity, assuming that there are no transaction costs.

\subsection{Harmonic analysis on semigroups}
We start by a brief introduction on harmonic analysis on
semigroups, for a complete treatment see \cite{Berg84b} and the
references therein. Unless otherwise specified, all measures are
supposed to be positive.

\begin{definition}
A function $\rho:\mathbb{S}\rightarrow\reals$ is called a
\textit{semicharacter} iff it satisfies $\rho(st)=\rho(s)\rho(t)$
for all $s,t \in \mathbb{S}$ and $\rho(1_{\mathbb{S}})=1$.
\end{definition}

In \cite{Berg84b} an involution operation is defined on the
semigroup $(\mathbb{S},\cdot)$, here and in the rest of the paper
we suppose that involution to be the identity, which means in
particular that we take all semicharacters to be real valued. The
dual semigroup of $\mathbb{S}$, i.e. the set of semicharacters on
$\mathbb{S}$ is called $\mathbb{S}^*$. In this context, we call a
function $f:\mathbb{S}\rightarrow\reals$ a \textit{moment function
on $\mathbb{S}$} iff $f(1_{\mathbb{S}})=1$ and $f$ can be
represented as:
\BEQ
f(s)=\int_{\mathbb{S}^*}{\rho(s)d \nu(\rho)},\quad\mbox{for all
}s\in\mathbb{S}, \label{repres-moment}
\EEQ
where $\nu$ is a Radon measure on $\mathbb{S}^*$.

When $\mathbb{S}$ is the semigroup defined in
(\ref{ss-market-structure}) as an enlargement of the semigroup of
monomials on $\reals^n$, its dual $\mathbb{S}^*$ is the set of
applications $\rho_x:\mathbb{S}\rightarrow\reals$ such that
$\rho_x(s)=s(x)$ for all $s\in\mathbb{S}$ and all $x\in\reals^n$.
The measure $\nu$ is then assimilated to a probability measure on
$\reals^n$ and the representation above becomes:
\BEQ
f(s)=\textstyle \Expect_{\nu}\left[s(x)\right],\quad\mbox{for all
}s\in\mathbb{S}. \label{repres-prob}
\EEQ
Our objective below is to find tractable conditions for a set of
prices $p_0,\ldots,p_{n+m}$ to be represented as
$\Expect_{\nu}\left[|w_i^{T}x-K_i|\right]=p_i$ for
$i=0,\ldots,n+m$ and some positive measure $\nu$.

\subsection{The compact case}
In this section we assume the asset distribution has a compact
support $K$. We treat the compact case independently as it is
rather simple yet captures many of the key features of the general
result. We begin by a few definitions along the lines of
\cite{Berg84a} and \cite{Berg84b}. An \textit{absolute value} on
$\mathbb{S}$ is a function $|\cdot|:\mathbb{S} \rightarrow
\reals_+$ satisfying
\[
|s^2| \leq |s|^2,\quad \mbox{for all } s\in\mathbb{S}
\]
and
\[
|1_{\mathbb{S}}|\geq 1.
\]
A function $f:\mathbb{S} \rightarrow \reals$ is said to be bounded
with respect to an absolute value $|\cdot|$ iff there exists some
$M>0$ such that
\[
|f(s)|\leq M |s|,\quad \mbox{for all } s\in\mathbb{S}.
\]
Furthermore, $f$ is called \textit{exponentially bounded} iff $f$
is bounded with respect to some absolute value. Remark that if the
measure $\nu$ in (\ref{repres-moment}) has its support contained
in the compact $K$ then the moment function
$f(s)=\int_{\mathbb{S}^*}{\rho(s)d \nu(\rho)}$ is bounded with
respect to the following absolute value:
\[
|s|_K=\sup_{\rho \in K}{\rho(s)}
\]
for $s \in \mathbb{S}$.

\begin{definition}
A function $f:\mathbb{S} \rightarrow \reals$ is called positive
semidefinite iff for all finite families $\{s_i\}$ of elements of
$\mathbb{S}$, the matrix with coefficients $f(s_i s_j)$ is
positive semidefinite.
\end{definition}

We remark that moment functions are necessarily positive
semidefinite. Necessary and sufficient conditions for the
existence of a measure $\nu$ in (\ref{repres-prob}) were derived
in \cite{henk90}, they were however numerically intractable. Here,
based on the results in \cite{Berg84b}, \cite{Puti99} and
\cite{Roma03}, we look for exploitable conditions for
representation (\ref{repres-prob}) to hold.

Let $\alpha$ be an absolute value, the central result in \cite[Th.
2.6]{Berg84b} states that the set of $\alpha$-bounded positive
semidefinite functions $f:\mathbb{S}\rightarrow \reals$ such that
$f(1_{\mathbb{S}})=1$ is a Bauer simplex whose extreme points are
given by the set of $\alpha$-bounded semicharacters. Hence a
function $f$ is positive semidefinite and exponentially bounded if
and only if it can be represented as
$f(s)=\int_{\mathbb{S}^*}{\rho d\nu(\rho)}$ with the support of
$\nu$ included in some compact subset of $\mathbb{S}^*$.

Based on these results, we derive below a set of tractable
necessary and sufficient conditions allowing a function $f$ to be
represented as in (\ref{repres-prob}). For $s,u$ in $\mathbb{S}$,
we note $E_s$ the shift operator such that for $f:\mathbb{S}
\rightarrow \reals$, we have $E_s(f(u))=f(su)$ and we let
$\mathcal{E}$ be the commutative algebra generated by the shift
operators on $\mathbb{S}$. Finally, we let $\beta=\sup_{x\in
K}\{\sum_{i=0}^{n+m}{e_i(x)}\}$.

\begin{theorem}
Suppose the asset distribution has compact support $K$ and
$\mathbb{S}$ is the payoff semigroup defined in
(\ref{ss-market-structure}), with $\beta$ is defined as above. A
function $f(s):\mathbb{S} \rightarrow \reals$ can be represented
as
\BEQ
f(s)=\textstyle \Expect_{\nu}[s(x)],\quad\mbox{for all
}s\in\mathbb{S}, \label{th:repres-prob}
\EEQ
for some measure $\nu$ on $K$, and satisfies the price constraints
in (\ref{eq:constraints}) if and only if:
\begin{enumerate}
\item[(i)] $f$ is positive semidefinite,
\item[(ii)] $E_{e_i}f$ is positive semidefinite for $i=0,\ldots,n+m,$
\item[(iii)] $\left(\beta I - \sum_{i=0}^{n+m}{E_{e_i}}\right)f$ is positive
semidefinite,
\item[(iv)] $f(e_{i})=p_i$ for $i=1,\ldots,n+m.$
\end{enumerate}
Furthermore, for each function $f$ satisfying conditions (i) to
(iv), the measure $\nu$ in representation (\ref{th:repres-prob})
is unique.
\end{theorem}
\begin{proof}
The family of shift operators
$\tau=\{\{E_{e_i}\}_{i=0,\ldots,n+m},\left(\beta I -
\sum_{i=0}^{n+m}{E_{e_i}}\right)\}\subset\mathcal{E}$ is such that
$I-T\in\mathrm{span}^+\tau$ for each $T\in\tau$ and
$\mathrm{span}\;\tau=\mathcal{E}$, hence $\tau$ is linearly
admissible in the sense of \cite[Corollary 2.5]{Berg84a} or
\cite{Mase77}, which states that (ii) and (iii) are equivalent to
$f$ being $\tau$-positive. Then, \cite[Th. 2.1]{Mase77} means that
$f$ is $\tau$-positive if and only if there is a measure $\nu$
such that $f(s)=\int_{\mathbb{S}^*}{\rho(s)d \nu(\rho)}$, whose
support is a compact subset of the $\tau$-positive semicharacters.
This means in particular that for a semicharacter $\rho_x \in
\mathrm{supp}(\nu)$ we must have $\rho_x(e_i)\geq0$, for
$i=1,\ldots,n$ hence $x\geq0$. The set of $\tau$-positive
semicharacters is then included in the nonnegative orthant and
includes both the simplex $\{x\geq0:\|x\|_1\leq\beta\}$ and K,
hence $f$ being $\tau$-positive is equivalent to $f$ admitting a
representation of the form $f(s)=\Expect_{\nu}\left[s(x)\right]$,
for all $s\in\mathbb{S}$ with $\nu$ having a compact support
$K\subset\reals_+^n$.
\end{proof}

\subsection{The unbounded case}
\label{ss-unbounded-case} The conditions derived in the last part
do not describe all possible arbitrage free prices as they cannot
account for unbounded asset distributions. Here, we use results
from \cite{Puti99} and \cite{Roma03} to derive intrinsic
characterizations of viable multivariate straddle prices.

We note $\mathcal{A}(\mathbb{S})$ the $\reals$-algebra generated
by the functions $\chi_s:\mathbb{S}^*\rightarrow\reals$ such that
$\chi_s(\rho)=\rho(s)$ for all $s\in\mathbb{S}$. By construction,
$\chi_s(\rho)=E_s\rho(1_{\mathbb{S}})$, and for a polynomial
$p\in\mathcal{A}(\mathbb{S})$ with $p=\sum_k q_k \chi_{g_k}$ and
for $\rho\in\mathbb{S}^*$ we have $p\rho(s)=\sum_k q_k \rho(s
g_k)$ for all $s\in\mathbb{S}$. When $\mathbb{S}$ is the payoff
semigroup defined in (\ref{ss-market-structure}), we naturally
have $\chi_{s}(\rho_x)=s(x)$, for all $x\in\reals^n$, $s\in
\mathbb{S}$ and $\rho \in \mathbb{S}^*$.

We now note $\mathcal{A}_\theta(\mathbb{S})$ the $\reals$-algebra
generated by $\mathcal{A}(\mathbb{S})$ and $\theta$ where
\BEQ
\theta(\rho)=\left(1+\sum_{i=0}^{m+n}{\chi_{e_i^2}(\rho)}\right)^{-1},\quad\mbox{for
all }\rho\in\mathbb{S}^*, \label{f-theta}
\EEQ
we also note $\mathcal{A}(\mathbb{S},y)$ the algebra generated by
$\mathcal{A}(\mathbb{S})$ and $\reals[y]$. We first simplify the
equality constraints on $2^n$ variables in \cite[Th. A]{Roma03} to
recover an additive formulation as in \cite{Puti99}. We begin by
proving the following lemma.

\begin{lemma}
\label{lemma-kernel} The kernel of the algebra homomorphism
$\Phi$:
\BEQ
\BA{lcl}
\mathcal{A}(\mathbb{S},y) & \rightarrow &
\mathcal{A}_\theta(\mathbb{S}) \\
p(\rho,y) & \mapsto & \Phi p=p(\rho,\theta(\rho))
\EA
\EEQ
is the ideal generated by $\sigma\in\mathcal{A}(\mathbb{S},y)$
such that
$\sigma(\rho,y)=y(1+\sum_{i=1}^{m+n+1}{\chi_{e_i^2}(\rho)})-1$.
\end{lemma}
\begin{proof}
We adapt the proof of \cite[lemma 2.3]{Puti99} and let $p \in
\mathcal{A}(\mathbb{S},y)$ be such that $p(\rho,\theta(\rho))=0$,
we write $p(\rho,y)=\sum_{k}{q_k(\rho)}y^k$ with $q_k \in
\mathcal{A}(\mathbb{S})$. We have:
\[
\BA{ll}
p(\rho,y) & =p(\rho,y)-p(\rho,\theta(\rho))= \sum_{k>0}{q_k(\rho)}(y^k-(\theta(\rho))^k) \\
          & =(y-(\theta(\rho)))l(\rho,y,(\theta(\rho)),
\EA
\]
where $l$ is a polynomial. Let $\kappa=max\{k:q_k\neq0\}$ and
\[
\tau(\rho)=\left(1+\sum_{i=1}^{m+n+1}{\chi_{e_i^2}(\rho)}\right)^\kappa,\quad\mbox{for
all }\rho\in\mathbb{S}^*,
\]
we then have
\BEQ
\tau(\rho)p(\rho,y)=\sigma(\rho,y)r(\rho,y),
\label{proof-tau-identity}
\EEQ
with $r(\rho,y) \in \mathcal{A}(\mathbb{S},y)$. The case
$\kappa=0$ is trivial hence we can assume $\kappa\neq0$. Using the
fact that the polynomials $\tau(z)$ and $\sigma(z)$ have no common
zeroes in $\mathbf{C}^{m+n+2}$[z], Hilbert's
\textit{Nullstellensatz} (see \cite{boch98} for example) states
that there must be
$\tilde{\tau},\tilde{\sigma}\in\mathbf{C}^{m+n+2}[z]$ such that
\[
\tau\tilde{\tau}+\sigma\tilde{\sigma}=1.
\]
Multiplying this last identity by $p$ yields, together with
(\ref{proof-tau-identity}):
\[
p=\sigma(r\tilde{\tau}+p\tilde{\sigma})
\]
hence the desired result.
\end{proof}

The next proposition is adapted from the dimensional extension
method in \cite[Th. 2.5]{Puti99} and \cite[Th. 4]{Roma03}, to
replace the exponential number of equality constraints in
\cite[Th. A]{Roma03} with an additive formulation as in
\cite{Puti99}. The function $\theta(\rho)$ is defined as in
(\ref{f-theta}) and $\mathcal{A}_\theta(\mathbb{S})$ is the
$\reals$-algebra generated by $\mathcal{A}(\mathbb{S})$ and
$\theta$.

\begin{proposition}
\label{prop-measure} With $\mathbb{S}$ being the payoff semigroup
defined in (\ref{ss-market-structure}), let $\Lambda$ be a
positive semidefinite linear form on
$\mathcal{A}_\theta(\mathbb{S})$ such that $\Lambda\left(x_i r^2
\right) \geq 0$ for all $r \in \mathcal{A}_\theta(\mathbb{S})$ and
$i=1,\ldots,n$, then $\Lambda$ has a unique representing measure
$\nu$ with support in $\reals^n_+$ and
$\mathcal{A}_\theta(\mathbb{S})$ is dense in $L^2(\nu)$.
\end{proposition}
\begin{proof}
We recall that the linear form $\Lambda$ is positive semidefinite
iff $\Lambda(r^2)\geq$, for all
$r\in\mathcal{A}_\theta(\mathbb{S})$. As in \cite{Roma03}, we
define a bilinear form on $r\in\mathcal{A}_\theta(\mathbb{S})$ by:
\[
\langle r_1,r_2 \rangle:=\Lambda(r_1 r_2),\quad \mbox{for all }
r_1,r_2 \in \mathcal{A}_\theta(\mathbb{S})
\]
We let $\mathcal{N}$ be the set $\{r \in
\mathcal{A}_\theta(\mathbb{S}):\Lambda(r^2)=0\}$. The bilinear
form above then defines a scalar product on
$\mathcal{A}_\theta(\mathbb{S})/\mathcal{N}$, and we note
$\mathcal{H}$ the completion of this space. We define in
$\mathcal{H}$ the operators:
\[
T_i(r+\mathcal{N})=\chi_{e_i}r+\mathcal{N},\quad \mbox{for all }r
\in \mathcal{A}_\theta(\mathbb{S})/\mathcal{N} \mbox{ and }
i=0,\ldots,n+m,
\]
which are symmetric and densely defined in $\mathcal{H}$. We also
define the operator $\left(D(B),B\right)$ by:
\[
D(B)=\mathcal{A}_\theta(\mathbb{S})/\mathcal{N} \mbox{ and }
B=\sum_{i=0}^{m+n}{T_i^2}.
\]
The operator $B$ is positive as a sum of squares of operators and,
by construction, the domain $D(B)$ is dense in $\mathcal{H}$ and
invariant by $B$. Let $\tau=\sum_{i=0}^{m+n}{\chi_{e_i^2}(\rho)}$
and $r \in \mathcal{A}_\theta(\mathbb{S})/\mathcal{N}$, then
$u=r\theta$ is such that $(1+\tau)u=r$, hence the operator $I+B$
is bijective on $D(B)$. This means that $B$ satisfies the
hypothesis of \cite[Lemma 2.2]{Puti99} and is essentially
self-adjoint. \cite[Prop. 1]{Roma03} then implies that the
operators $T_i$ for $i=0,\ldots,n+m$ are essentially normal and
that their canonical closures commute, meaning that there exists a
common spectral measure $H$ for the operators $\bar T_i$ for
$i=0,\ldots,n+m$. With $T=(T_i)_{i=0,\ldots,n+m}$ and $r \in
\mathcal{A}_\theta(\mathbb{S})/\mathcal{N}$, we define the
operator $r(T)$ by:
\BEQ
\BA{lcl}
\mathcal{A}_\theta(\mathbb{S})/\mathcal{N} & \rightarrow &
\mathcal{A}_\theta(\mathbb{S})/\mathcal{N} \\
w+\mathcal{N} & \mapsto & r(T)(w+\mathcal{N})=rw+\mathcal{N}.
\EA
\EEQ
With $\gamma(x)=\left(\sum_{i=0}^{n+m}{x_i^2}\right)^{-1}$, there
is an element $q$ of $\reals_\gamma[x]$, the $\reals$-algebra
generated by $\reals[x]$ and $\gamma(x)$ such that
$r(\rho)=q\left((\chi_{e_i})_{i=0,\ldots,n+m}(\rho),\theta(\rho)\right)$
for all $\rho \in \mathbb{S}^*$. We then have:
\[
\Lambda(r)=\langle r(T)1,1 \rangle=\langle q(\bar T)1,1
\rangle=\int_{\mathbf{R}^{n+m+1}}q(x)dH_{1+\mathcal{N},1+\mathcal{N}}(x),
\]
The homomorphism $f$:
\BEQ
\BA{lcl}
\reals_\gamma[x] & \rightarrow &
\mathcal{A}_\theta(\mathbb{S}) \\
p(x) & \mapsto &
f(p)=p((\chi_{e_i}(\rho))_{i=0,\ldots,n+m},\theta(\rho))
\EA
\EEQ
satisfies the hypothesis of \cite[Lemma 2]{Roma03} hence there is
a (positive) Radon measure $\nu$ on such that:
\[
\Lambda(r)=\int_{\mathbb{S}^*}{r(\rho)d \nu(\rho)},
\]
which, if $\mathbb{S}$ is defined as in section
(\ref{ss-market-structure}), is also:
\[
\Lambda(r)=\int_{\mathbf{R}^n}{r(x)d \nu(x)}.
\]
Uniqueness and density follow from the argument in \cite{Roma03}.
Now, because the operators $T_i$ for $i=1,\ldots,n$ are
essentially self-adjoint with $\Lambda\left(x_i r^2 \right) \geq
0$ for $r \in \mathcal{A}_\theta(\mathbb{S})$ and $i=1,\ldots,n$,
we know that the $T_i$ are positive for all $i$. The spectral
measure $F_i$ of $T_i$ is given by $F_i(X)=H\left(\bar
T_i^{-1}(X)\right)$ for all Borel sets $X \subset\reals$ and $F_i$
must be concentrated in $\reals_+$ for all $i=1,\ldots,n$ hence
the spectral measure $H$ of $\bar T$ is concentrated in
$\reals^n_+$ and so is the representing measure $\nu$.
\end{proof}

We can now formulate a general moment theorem that describes all
the price systems that admit a representation as in
(\ref{repres-prob}).

\begin{theorem}
\label{th:repres-prob-unbounded} Let $\mathbb{S}$ be defined as in
(\ref{ss-market-structure}). A sequence $f(s):\mathbb{S}
\rightarrow \reals$ is a moment sequence and can be represented as
in (\ref{repres-prob}):
\[
f(s)=\textstyle \Expect_{\nu}[s(x)],\quad\mbox{for all
}s\in\mathbb{S},
\]
for some measure $\nu$ with support in $\reals_+^n$, if and only
if there is a sequence $p(s,k):\left( \mathbb{S},\mathbf{N}
\right) \rightarrow \reals$ such that:
\begin{enumerate}
\item[(i)] $p(s,0)=f(s)$ for all $s \in \mathbb{S}$,
\item[(ii)] $p(s,k)$ is positive semidefinite on
$\left(\mathbb{S},\mathbf{N}\right)$,
\item[(iii)] $p(e_i s,k)$ is positive semidefinite on $\left(\mathbb{S}
,\mathbf{N}\right)$ for $i=1,\ldots,n$,
\item[(iv)] $p(s,k)=p(s,k+1)-\sum_{i=0}^{n+m}{p(e_i^2s,k+1)}$ for
all $(s,k) \in (\mathbb{S},\mathbf{N})$.
\end{enumerate}
Furthermore, the representing measure for sequence $f$ is unique
if and only if the sequence $p$ is unique.
\end{theorem}
\begin{proof}
First we show that conditions (i)-(iv) are necessary. With
$\mathbb{S}$, the payoff semigroup defined in
(\ref{ss-market-structure}), we recall that $\mathbb{S}^*$ can be
identified with $\reals_+^n$, hence $\chi_{s}(\rho_x)=s(x)$, for
all $x\in\reals_+^n$, $s\in \mathbb{S}$ and $\rho \in
\mathbb{S}^*$. Suppose that $f$ can be represented as:
\[
f(s)=\int_{\mathbf{R}_+^n}{s(x)d \nu(x)},\quad \mbox{for all }s
\in \mathbb{S},
\]
we let
\[
p(s,k)=\int_{\mathbf{R}^n_+}{s(x)
\left(1+\sum_{i=0}^{m+n}{e_i^2}(x)\right)^{-k} d \nu(x)},\quad
\mbox{for all }(s,k) \in (\mathbb{S},\mathbf{N}),
\]
which satisfies (i) and (iv) by construction, $p(s,k)$ is then a
moment sequence on the product semigroup
$\left((\mathbb{S},\cdot)\times(\mathbf{N},+)\right)$ and as such
must be positive semidefinite, hence condition (ii). Then, because
for $i=1,..,n$ we have
\[
p(e_is,k)=\int_{\mathbf{R}^n_+}{s(x)
\left(1+\sum_{i=0}^{m+n}{e_i^2}(x)\right)^{-k} e_i(x) d
\nu(x)},\quad \mbox{for all }(s,k) \in (\mathbb{S},\mathbf{N}),
\]
we know that $p(e_is,k)$ is a moment sequence for the measure
$e_i(x)d \nu$, hence condition (iii).

Conversely, let's assume that we are given a sequence $p(s,k)$
satisfying (i)-(iv). We let $\mathcal{A}_\theta(\mathbb{S})$ and
$\mathcal{A}(\mathbb{S},y)$ be the $\reals$-algebras described at
the beginning of the section. We define a linear function
$\Lambda$ on $\mathcal{A}(\mathbb{S},y)$ by:
\[
L\left(
\sum_{j,k}{a_j\chi_{s_j}y^k}\right)=\sum_{j,k}{a_jp(s_j,k)}
\]
and as in lemma \ref{lemma-kernel}, we can define the following
algebra homomorphism $\Phi$:
\BEQ
\BA{lcl}
\mathcal{A}(\mathbb{S},y) & \rightarrow &
\mathcal{A}_\theta(\mathbb{S}) \\
p(\rho,y) & \mapsto & \Phi p=p(\rho,\theta(\rho))
\EA
\EEQ
whose kernel $\mathcal{N}$ has been computed in lemma
\ref{lemma-kernel}, and $\mathcal{A}_\theta(\mathbb{S})$ is
isomorphic to the quotient
$\mathcal{A}(\mathbb{S},y)/\mathcal{N}$. Condition (iv) implies
that $L(\mathcal{N})=0$ and we can then define a linear form
$\Lambda$ on $\mathcal{A}_\theta(\mathbb{S})$ by:
\[
\Lambda(r)=L(q),\quad \mbox{where
}r(\rho)=q(\rho,\theta(\rho)),\quad \mbox{for all }\rho \in
\mathbb{S}^*,
\]
with $r \in \mathcal{A}_\theta(\mathbb{S})$ and $q \in
\mathcal{A}(\mathbb{S},y)$. Because of (i)-(iv), the form
$\Lambda$ satisfies the hypothesis of proposition
\ref{prop-measure} and has a unique representing measure $\nu$.
\end{proof}

\section{Price bounds and static hedging}
In this section, we show how the duality between the existence of
a pricing measure and that of a replicating portfolio transposes
into the moment framework described in the previous section. In
particular, we detail how an optimal static super/sub-replicating
portfolio can be constructed using the solution to the dual of to
the moment problem in (\ref{eq:constraints}). In particular, in a
result that is consistent with the dynamic framework (see
\cite{Avel95}), the replicating portfolio only involves options in
the data set and no other option is needed to "complete the grid".

\subsection{Price bounds via semidefinite programming}
Here, we show how one can compute bounds on the solution of
problem (\ref{eq:constraints}) using a subset of the moment
conditions imposed by theorem \ref{th:repres-prob-unbounded}.
These conditions cast (\ref{repres-prob}) as a semidefinite
program (see \cite{Nest94} or \cite{Vand96}), which can then be
solved efficiently using solvers such as SEDUMI by \cite{stur99}.

\subsubsection{Asset distributions with compact support}
As before, we note $\mathcal{A}(\mathbb{S})$ the $\reals$-algebra
generated by the functions $\chi_s:\mathbb{S}^*\rightarrow\reals$
such that $\chi_{s}(\rho)=\rho(s)$ for all $s\in\mathbb{S}$ and
$\rho\in\mathbb{S}^*$. For a polynomial
$p\in\mathcal{A}(\mathbb{S})$ with $p=\sum_i q_i \chi_{g_i}$ where
$g_i \in \mathbb{S}$, and for $s \in \mathbb{S}$ we set
\[
p\rho(s)=\sum_i q_i \rho(sg_i).
\]
With $\mathbb{S}$ the payoff semigroup defined in
(\ref{ss-market-structure}), we recall that $\mathbb{S}^*$ can be
identified with $\reals^n$, hence $\chi_{s}(\rho_x)=s(x)$, for all
$x\in\reals^n$, $s\in \mathbb{S}$ and $\rho_x \in \mathbb{S}^*$.
This means that $p\in\mathcal{A}(\mathbb{S})$ can be rewritten
\[
p(x)=\sum_i q_i s(x)g_i(x), \quad \mbox{for all }x\in\reals_+^n.
\]
We now recall the construction of moment matrices as in
\cite{Curt00} and \cite{Lass01}. We adopt the following multiindex
notation for monomials in $\mathcal{A}(\mathbb{S})$:
\[
e^{\alpha}(x):=e_0^{\alpha_0}(x)e_1^{\alpha_1}(x)\cdots
e_{m+n}^{\alpha_{m+n}}(x),
\]
and we let
\BEQ
y_e=(1,e_0,\ldots,e_{m+n},e_0^2,e_0e_1,\ldots,
e_0^d,\ldots,e_{m+n}^d) \label{def:moment-vector-compact}
\EEQ
be the vector of all monomials in $\mathcal{A}(\mathbb{S})$, up to
degree $d$, listed in graded lexicographic order. We note $s(d)$
the size of the vector $y_e$. Let $y\in \reals^{s(2d)}$ be the
vector of moments (indexed as in $y_e$) of some probability
measure $\nu$ with support in $\reals^n_+$, we note $M_{d}(y)\in
\reals^{s(d)\times s(d)}$, the symmetric matrix:
\[
M_{d}(y)_{i,j}=\int_{\mathbf{R}^n_+}\left(y_{e}\right)_i (x)
\left(y_{e}\right)_j (x) d\nu(x),\quad \text{for }i,j=1,...,s(d)
\]
In the rest of the paper, we will always implicitly assume that
$y_1=1$. With $\beta (i)$ the exponent of the monomial
$\left(y_{\chi}\right)_i$ and conversely, $i(\beta)$ the index of
the monomial $e^{\beta}$ in $y_{e}$. We notice that for a given
moment vector $y\in \reals^{s(d)}$ ordered as in
(\ref{def:moment-vector}), the first row and columns of the matrix
$M_{d}(y)$ are then equal to $y$. The rest of the matrix is then
constructed according to:
\[
M_{d}(y)_{i,j}=y_{i(\alpha +\beta)}\text{ if
}M_{d}(y)_{i,1}=y_{i(\alpha)}\text{ and
}M_{d}(y)_{1,j}=y_{i(\beta)}.
\]

Similarly, let $g\in \mathcal{A}(\mathbb{S})$, we derive the
moment matrix for the measure $g(x)d\nu$ on $\reals_+^n$ (called
the localizing matrix in \cite{Curt00}), noted $M_{d}(gy)\in
\symm^{s(d)}$, from the matrix of moments $M_{d}(y)$ by:
\[
M_{d}(gy)_{i,j}=\int_{\mathbf{R}^n_+}\left( y_{e}\right)_{i}(x)
\left(y_{e}\right)_{j}(x) g(x) d\nu(x)
\]
for $i,j=1,...,s(d)$. The coefficients of the matrix $M_{m}(gy)$
are then given by:
\BEQ
M_{d}(gy)_{i,j}=\sum_{\alpha}g_{\alpha}y_{i(\beta
(i)+\beta(j)+\alpha)} \label{def:localizing-matrix-compact}
\EEQ

We can then form a semidefinite program to compute a lower bound
on the optimal solution to (\ref{eq:constraints}) using a subset
of the moment constraints in theorem \ref{th:repres-prob}, taking
only monomials and moments in $y$ up to a certain degree.
\begin{corollary}
\label{th:compact-program} Let $N$ be a positive integer and
$y\in\reals^{s(2N)}$, a lower bound on the optimal value of:
\[
\BA{ll}
\mbox{minimize} & {p_0:=\Expect}_{\nu }[e_0(x)] \\
\mbox{subject to} & {\Expect}_{\nu }[e_i(x)]=p_i, \quad
i=1,\ldots,n+m,
\EA
\]
can be computed as the solution of the following semidefinite
program:
\BEQ
\BA{ll}
\mbox{minimize} & y_2 \\
\mbox{subject to} & M_N(y) \succeq 0 \\
 & M_N(e_jy) \succeq 0 , \quad \mbox{for }j=1,\ldots,n,\\
 & M_N\left((\beta-\sum_{k=0}^{n+m}{e_k})y)\right) \succeq 0 \\
 & y_{(j+2)}=p_j, \quad \mbox{for }j=1,\ldots,n+m \mbox{ and }
 s\in \mathbb{S}\\
\EA
\label{eq:sdp-bound-compact}
\EEQ
where $s$ is such that $i(s)\leq s(2N)$. The optimal value of
(\ref{eq:sdp-bound-compact}) converges to the optimal value of the
original program as $N \rightarrow \infty$.
\end{corollary}

\subsubsection{Unbounded distributions}
Here we work on the product semigroup
$(\mathbb{S},\cdot)\times(\mathbf{N},+)$. Its dual is the set of
functions $\rho_x:(\mathbb{S},\mathbf{N})\rightarrow\reals$ such
that $\rho_x((s,k))=s(x)x^k$ for all $s\in\mathbb{S}$,
$k\in\mathbf{N}$ and $x\in\reals^n$. As before, we note
$\mathcal{A}(\mathbb{S},\mathbf{N})$ the $\reals$-algebra
generated by the functions
$\chi_s:(\mathbb{S},\mathbf{N})^*\rightarrow\reals$ such that
$\chi_{(s,k)}(\rho)=\rho((s,k))$ for all $s\in\mathbb{S}$ and
$\rho\in(\mathbb{S},\mathbf{N})^*$. With $\mathbb{S}$, the payoff
semigroup defined in (\ref{ss-market-structure}), here
$(\mathbb{S},\mathbf{N})^*$ can again be identified with
$\reals^n$, hence $\chi_{(s,k)}(\rho_x)=s(x)x^k$, for all
$x\in\reals^n$, $(s,k)\in (\mathbb{S},\mathbf{N})$ and $\rho_x \in
(\mathbb{S},\mathbf{N})^*$. By construction, we have
\[
(\chi_{(s,k)})^2=\chi_{(s^2,2k)},\quad\mbox{for all }(s,k)
\in(\mathbb{S},\mathbf{N}),
\]
and for a polynomial $p\in\mathcal{A}(\mathbb{S},\mathbf{N})$ with
$p=\sum_i q_i \chi_{g_i}x^{k_i}$ where $(g_i,k_i) \in
(\mathbb{S},\mathbf{N})$, and for $(s,l) \in
(\mathbb{S},\mathbf{N})$ we set
\[
p((s,l))(x)=\sum_i q_i s(x)g_i(x)x^{k_i+l},
\]
for all $x\in\reals^n$. We adopt here the multiindex notation for
monomials in $\mathcal{A}(\mathbb{S},\mathbf{N})$:
\[
e^{\alpha}:={(e_0,0)}^{\alpha_0}{(e_1,0)}^{\alpha_1}\cdots
{(e_{m+n},0)}^{\alpha_{m+n}}{(1,1)}^{\alpha_{m+n+1}}.
\]
We then let
\BEQ
y_e=(1,{(e_0,0)},\ldots,{(e_{m+n},0)},{(1,1)},
{(e_0,0)}^2,{(e_0,0)}{(e_1,0)},\ldots,
{(e_0,0)}^d,\ldots,{(1,1)^d}) \label{def:moment-vector}
\EEQ
be the vector of all monomials in
$\mathcal{A}(\mathbb{S},\mathbf{N})$, up to degree d, listed in
graded lexicographic order. We note $s(d)$ the size of the vector
$y_e$. The matrices $M_d(y)$ and $M_d(gy)$ are defined as in the
compact case above.

We can again form a semidefinite program, this time using a subset
of the moment constraints in theorem
\ref{th:repres-prob-unbounded}, taking only moments up to a
certain degree.
\begin{corollary}
Let $N$ be a positive integer and $y\in\reals^{s(2N)}$, a lower
bound on the optimal value of:
\[
\BA{ll}
\mbox{minimize} & {p_0:=\Expect}_{\nu }[e_0(x)] \\
\mbox{subject to} & {\Expect}_{\nu }[e_i(x)]=p_i, \quad
i=1,\ldots,n+m,
\EA
\]
can be computed as the solution of the following semidefinite
program:
\BEQ
\BA{ll}
\mbox{minimize} & y_2 \\
\mbox{subject to} & M_N(y) \succeq 0 \\
 & M_N((e_j,0)y) \succeq 0 ,\quad \mbox{for }j=1,\ldots,n,\\
 & y_{i(s,k)}=y_{i(s,k+1)}-\sum_{i=0}^{n+m}{y_{i(e_i^2s,k+1)}} \\
 & y_{(j+1)}=p_j, \quad \mbox{for }j=1,\ldots,n+m \mbox{ and }
 (s,k)\in (\mathbb{S},\mathbf{N})\\
\EA
\label{eq:sdp-bound}
\EEQ
where $(s,k)$ are taken such that $i(s,k)\leq s(2N)$. The optimal
value of (\ref{eq:sdp-bound}) converges to the optimal value of
the original program as $N \rightarrow \infty$.
\end{corollary}

\subsection{Static hedging portfolios and sums of squares}
We let here $\Sigma\subset\mathcal{A}(\mathbb{S})$ be the set of
polynomials that are sums of squares of polynomials in
$\mathcal{A}(\mathbb{S})$, and $\mathcal{P}$ the set of positive
semidefinite sequences on $\mathbb{S}$. The central argument of
this paper is to replace the conic duality between probability
measures and positive portfolios:
\[
p(x)\geq0 \Leftrightarrow \int p(x)d\nu\geq0,\quad\mbox{for all
measures }\nu,
\]
by the conic duality between positive semidefinite sequences
$\mathcal{P}$ and sums of squares polynomials $\Sigma$:
\[
\langle f,p \rangle \geq 0 \mbox{ for all $p\in\Sigma$ iff }f\in
\mathcal{P}.
\]
for $p\in\mathcal{A}(\mathbb{S})$ with $p=\sum_i q_i \chi_{s_i}$
and $f:\mathbb{S}\rightarrow \reals$ having defined $\langle f,p
\rangle=\sum_i q_i f(s_i)$. The previous section used positive
semidefinite sequences to characterize viable price sets, in this
section, we use sums of squares polynomials to characterize
super/sub-replicating portfolios.

From the initial price problem (\ref{eq:constraints}) written in
terms of straddles:
\BEQ
\BA{ll}
\mbox{minimize} & p_0:=\int_{\mathbf{R}^n_+}e_0(x)d\nu(x) \\
\mbox{subject to} & \int_{\mathbf{R}^n_+}e_i(x)d\nu(x)=p_i, \quad
i=1,\ldots,n+m,\\
 & \int_{\mathbf{R}^n_+}d\nu(x)=1,
\EA
\EEQ
in the variable $\nu$, a positive measure on $\reals_+^n$. We can
form the Lagrangian:
\[
L(\nu,\lambda)=\lambda_{n+m+1}+\sum_{i=1}^{n+m}{\lambda_i
p_i}+\int_{\mathbf{R}^n_+}\left(e_0(x)-\sum_{i=1}^{n+m}\lambda_{i}e_i(x)-\lambda_{n+m+1}\right)d\nu(x)
\]
with variables $\nu$ and $\lambda \in \reals^{n+m+1}$. We obtain
the classic dual as a portfolio replication problem:
\BEQ
\BA{ll}
\mbox{maximize} & \lambda_{n+m+1} + \sum_{i=1}^{n+m}{\lambda_i p_i} \\
\mbox{subject to} &
e_0(x)-\sum_{i=1}^{n+m}\lambda_{i}e_i(x)-\lambda_{n+m+1} \geq 0,
\quad \mbox{for all }x \in \reals_+^n.
\EA
\label{eq:classic-dual}
\EEQ
in the variable $\lambda \in \reals^{n+m+1}$.

Unfortunately, the problem formulations above are numerically
intractable except in certain particular cases (see \cite{Bert00}
and \cite{dasp03b}). On the other hand, as we have seen in the
previous section, the conditions of theorem \ref{th:repres-prob}
turn problem (\ref{eq:constraints}) into an infinite dimensional
semidefinite program which can be relaxed to produce tractable
bounds on the solution of (\ref{eq:constraints}). Here, we detail
the accompanying duality theory to exhibit a static hedging
portfolios corresponding to these bounds.

We can assume without loss of generality that the payoff functions
$\{e_i(x)\}_{i=0,\ldots,m+n}$, together with the cash
$1_{\mathbb{S}}$, are linearly independent. Then \cite[Proposition
6.1.8 and Theorem 6.1.10]{Berg84b} hold and we can form a dual to
the cone of positive semidefinite functions on $\mathbb{S}$ as
follows. For $p\in\mathcal{A}(\mathbb{S})$ with $p=\sum_i q_i
\chi_{s_i}$ and $f:\mathbb{S}\rightarrow \reals$ with:
\[
\langle f,p \rangle=\sum_i q_i f(s_i),
\]
\cite[Theorem 6.1.10]{Berg84b} states that $\Sigma$ is the polar
cone of $\mathcal{P}$ for the above bilinear form, in other words:
\[
\langle f,p \rangle \geq 0 \mbox{ for all $p\in\Sigma$ iff }f\in
\mathcal{P}.
\]
We can use this conic duality to compute a dual to program
(\ref{eq:sdp-bound-compact}). Considering the compact case for
simplicity, Corollary (\ref{th:compact-program}) states that the
initial pricing problem:
\[
\BA{ll}
\mbox{minimize} & {p_0:=\Expect}_{\nu }[e_0(x)] \\
\mbox{subject to} & {\Expect}_{\nu }[e_i(x)]=p_i, \quad
i=1,\ldots,n+m,
\EA
\]
is equivalent to the following (infinite) semidefinite program:
\[
\BA{ll}
\mbox{minimize} & y_2 \\
\mbox{subject to} & M(y) \succeq 0 \\
 & M(e_jy) \succeq 0 , \quad \mbox{for }j=0,\ldots,n+m\\
 & M\left((\beta-\sum_{k=0}^{n+m}{e_k})y)\right) \succeq 0 \\
 & y_{(j+2)}=p_j, \quad \mbox{for }j=1,\ldots,n+m\\
 & y_1=1.
\EA
\]
in the variable $y:\mathbb{S}\rightarrow \reals$. We can form the
Lagrangian:
\[
\BA{lll}
L(y,\lambda,q)&:=&y_2+(1-y_1)\lambda_{n+m+1}+\sum_{j=1}^{n+m}{(p_j-y_{(j+2)})\lambda_j} -\langle y,q_0 \rangle \\
&&- \sum_{j=0}^{n+m}{\langle e_jy,q_j \rangle}
- \langle (\beta-\sum_{k=0}^{n+m}{e_k})y,q_{n+1} \rangle\\
\EA
\]
or again:
\[
\BA{lll}
L(y,\lambda,q)&:=&y_2+(1-y_1)\lambda_{n+m+1}+\sum_{j=1}^{n+m}{(p_j-y_{(j+2)})\lambda_j} -\langle y,q_0 \rangle \\
&&- \sum_{j=0}^{n+m}{\langle y,e_jq_j \rangle}
- \langle y,(\beta-\sum_{k=0}^{n+m}{e_k})q_{n+1} \rangle\\
\EA
\]
in the variables $y:\mathbb{S}\rightarrow \reals$,
$\lambda\in\reals^{n+m+1}$ and $q_j\in\Sigma$ for
$j=0,\ldots,(n+1)$. We then get the dual as a portfolio problem:
\BEQ
\BA{ll}
\mbox{maximize} & \sum_{j=1}^{n+m}{p_j\lambda_j}+\lambda_{n+m+1} \\
\mbox{subject to} & e_0(x)-\sum_{j=1}^{n+m}{\lambda_je_j(x)}-\lambda_{n+m+1}\\
&=q_0(x)+\sum_{j=1}^{n+m}{q_j(x)e_j(x)}+(\beta-\sum_{k=0}^{n+m}{e_k(x)})q_{n+1}(x)
\EA
\label{eq:sdp-bound-compact-dual}
\EEQ
in the variables $\lambda\in\reals^{n+m+1}$ and $q_j\in\Sigma$ for
$j=0,\ldots,(n+1)$.

The key difference between this portfolio problem and the one in
(\ref{eq:classic-dual}) is that the (intractable) positivity
constraint
$e_0(x)-\sum_{i=1}^{n+m}\lambda_{i}e_i(x)-\lambda_{n+m+1} \geq 0$
in (\ref{eq:classic-dual}) is replaced by the tractable condition
that this portfolio be written as a combination of sums of squares
of polynomials in $\mathcal{A}(\mathbb{S})$. Such combinations can
be constructed directly from the dual solution to the semidefinite
program in (\ref{eq:sdp-bound-compact}), hence a numerical
solution to the program in (\ref{eq:sdp-bound-compact}) provides
both a price bound and an accompanying portfolio.

\bibliographystyle{amsalpha}
\providecommand{\bysame}{\leavevmode\hbox
to3em{\hrulefill}\thinspace}
\providecommand{\MR}{\relax\ifhmode\unskip\space\fi MR }
\providecommand{\MRhref}[2]{%
  \href{http://www.ams.org/mathscinet-getitem?mr=#1}{#2}
} \providecommand{\href}[2]{#2}

\end{document}